\documentclass{article}
\def \Z {{\mathbf {Z}}}

\def \N {{\mathbf {N}}}
\def \R {{\mathbf {R}}}
\def \B  {\cal B}

\def\uu{\bigsqcup}
\def\eps{\varepsilon}

\usepackage[T1]{fontenc} 
\usepackage[cp1251]{inputenc} 
\usepackage[russian]{babel}

\title{ Generic properties of ergodic automorphisms}
\author{V.V.~Ryzhikov}
\date{24.09.2024}

\begin{document}
\Large
\maketitle
\begin{abstract}
Typical properties of  measure space automorphisms with respect to the Halmos and Alpern-Tikhonov metrics are discussed.

\vspace{3mm}

UDK: 517.98      MSC: 28D05, 58F1
\end{abstract}

We consider the following metric spaces: the groups ${\bf Aut}$ and ${\bf Aut}_\infty$ with the Halmos metric and the families of mixing automorphisms $ {\bf Mix}$ and ${\bf Mix}_\infty$ with the Alpern-Tikhonov metric.
The group ${\bf Aut}$ of automorphisms of a probability space has been studied as a metric space since the 1940s \cite{Ha},\cite{Ro}. The set of mixing automorphisms $\bf Mix \subset Aut$, as Rokhlin showed, has the first Baire category with respect to the Halmos metric $\rho$. Alpern and Tikhonov proposed a stronger metric $r$, with respect to which the class $ {\bf Mix}$ turned out to be a complete space. This allows studying mixing actions by category methods. Ageev \cite{A05} gave a category proof of the existence of a weakly mixing automorphism with homogeneous spectrum of a given multiplicity
(for a simplified proof, see \cite{09}). Using the new metric, Tikhonov proved  the existence of  mixing automorphisms with such spectra  \cite{T11}. No explicit examples of mixing automorphisms with spectral multiplicity $m>2$ have been found.

Let ${\bf Aut}_\infty$ denote the group sigma-finite measure space  automorphisms. It is considered  with the Halmos metric $\rho_\infty$. The set of mixing automorphisms ${\bf Mix}_\infty$ will be  equipped with the Alpern-Tikhonov metric $r_\infty$.
The article discusses both known and some new facts about generic algebraic, spectral, approximation and entropy  properties of automorphisms. 
\section{The space ${\bf Aut}$, generic properties }
Recall the definition of Halmos metric $\rho$ on the group  all atomorphisms of the standard probability space $(X,\B,\mu)$:  
$$ \rho(S,T)=\sum_{i=0}^\infty 2^{-i}\left(\mu(SA_i\Delta TA_i)+\mu(S^{-1}A_i\Delta T^{-1}A_i)\right),$$
where $\{A_i\}$ is some fixed family dense in the algebra $\B$ of all $\mu$-measurable sets.
The space $({\bf Aut},\rho)$ is complete and separable. In ergodic theory we are interested in the question: which invariant properties  are generic? A property is generic (or called also typical) if the set of actions possessing it is generic.
The generic set, by definition, contains a dense $G_\delta$-set.
An automorphism of a  measure space  induces a unitary operator on $L_2$, whose spectral invariants become invariants of the automorphism. 

--- Simple spectrum, weak mixing, cyclic $(n,n+1)$-approximation,  the presence of all admissible weak limits (see, for example, \cite {Ka}, \cite {KS}, \cite {07}) are generic properties.

--- The symmetric tensor power $T^{\odot n}$ of the generic automorphism $T$ has  simple singular spectrum \cite{A99},
the local rank of the power $T^{\odot n}$ is $n!n^{-n}$ \cite{Ka},\cite{13}.

--- The generic automorphism has  singular spectrum because it is rigid ($T^{r_i}\to I$), but the Kushnirenko entropy of it with respect to growing sequences, for example, the $\{2^n\}$-entropy, is infinite \cite{21}. This contrasts with the properties of the horocyclic flow, which is multiple mixing, has  Lebesgue spectrum, and finite $\{2^n\}$-entropy \cite{Ku}. The generic automorphism has infinite $\alpha$-entropy for every sequence $\alpha$ such that $\alpha(n+1)- \alpha(n)\to+\infty$ \cite{21}.

--- The generic automorphism  is not isomorphic to its inverse \cite{J81} (explicit asymptotic invariants distinguishing them see in \cite{23}).

---  The theory of generic properties, starting with King's work \cite{Ki}, has acquired new methods and led to unexpected, subtle results. The generic automorphism has roots \cite{Ki}, moreover, it is included in the continuum of multidimensional non-isomorphic flows \cite{LR},\cite{SE},\cite{T06}. The centralizer of the  generic automorphism is studied in \cite{So}.

  --- Recall that  the  restriction of an action to a nontrivial invariant sigma-algebra  is called  the factor. The  generic automorphism has many factors, being simultaneously a compact extension  \cite{A00}, \cite{SE} and a relatively weakly mixing extension  \cite{GTW}. 

Will an automorphism that has a nongeneric  property be nongeneric?
This brings to mind the famous question: \it does an existing unicorn exist? \rm
Although we often say "general automorphism"{}, there is no such thing as a general automorphism. Generic properties  exist.

\vspace{2mm}
\bf Theorem 1.1. \it An automorphism that has all generic properties does not exist. \rm

\vspace{2mm}
Tikhonov noted that this statement  has an obvious proof. Indeed, a set that does not contain a given point is typical in our case, and the intersection of all such sets is empty. 

However we now  explain how to prove theorem 1.1  using  classical substantive invariants. At the dawn of ergodic theory, Halmos and Rokhlin published articles with contradictory titles: "In general, a measure-preserving transformation is mixing"{} \cite{Ha} and "A general transformation with an invariant measure is not  mixing"{} \cite{Ro}. In fact, there is no contradiction, since Halmos wrote \it mixing \rm instead of \it weak mixing. \rm For weakly mixing automorphism $T$ there is  some sequence ${n_i}$ such that 
$$ \mu(T^{n_i}A\cap B)\to \mu(A)\mu(B) \eqno (1)$$ for all measurable $A,B$.
An automorphism $T$ is (strongly) mixing if $$ \mu(T^{i}A\cap B)\to \mu(A)\mu(B), \ \ i\to \infty.$$ From the existence of  mixing automorphisms it follows that the  generic automorphism $T$ possesses a sequence $n_i\to\infty$ (it depends on $T$) such that (1) holds.
Thus,  generic $T$ is weakly mixing. Now we fix $T$ and its mixing sequence $n_i\to\infty$. For  generic $S$ there exists a subsequence ${n_{i(k)}}$ depending on $S$ such that  $S^{n_{i(k)}}\to I$, see \cite{21}. Thus, any generic $S$ is not mixing along  the sequence $n_i\to\infty$.

\vspace{2mm}
\bf Convergence of averages  for typical pairs $\bf S,T$. \rm Recently, in \cite{Au},
\cite{Ye},\cite{24aa} there appeared examples  of ergodic automorphisms
$S,T$ for which the averages  $ \sum_{n=1}^N  T^nf\,S^ng/N$
do not converge for some functions $f,g\in L_\infty$.

\vspace{2mm}
\bf Theorem 1.2. \it For a typical pair of automorphisms
$S,T$ and any $f,g\in L_\infty$  the averages  $$ \frac 1 N \sum_{n=1}^N  T^nf\,S^ng $$ converge in $L_2$.\rm

\vspace{2mm}
Proof. Let us consider a  set of pairs $S,T$ for which there exists a sequence
$n_i\to\infty$ which is  rigid for $S$ and mixing for $T$.
It is dense and  has type $G_\delta$. Such automorphisms $S,T$ are spectrally disjoint. This   implies the convergence   (see \cite{24aa})
$$ \frac 1 {N} \sum_{n=1}^{N} \int T^nf\,S^ng\, d\mu\ \to\
\int fd\mu \int gd\mu.\eqno (2) $$

Let $\nu$ be a measure on $X\times X\times X'\times X'$ ($X=X'$) obtained
as follows. For some sequence $N_i\to\infty$
for any real functions $f,g,f',g'\in L_\infty$ we have
$$ \frac 1 {N_i^2} \sum_{n,n'=1}^{N_i} \int T^nf\,S^ng\, d\mu\ \int T^{n'}f'\,S^{n'}g'\, d\mu\ \to\ \int f\otimes g \otimes f' \otimes g' d\nu.$$
From (2) it follows that the projection of the measure $\nu$ onto $X\otimes X$ (and $X'\otimes X'$) is $\mu\otimes \mu$.
But $\nu$ is invariant under the product $S\times T\times Id\times Id$,
and  $S\times T$ is ergodic for generic pairs $S,T$. Therefore,
$\nu=\mu\otimes \mu\otimes \mu\otimes \mu$. For $f$, $\int fd\mu=0$,
the scalar squares of  $A_N=\sum_{n=1}^{N}  T^nf\,S^ng/N$ tend to 0,
so $A_N$ converge to 0  in $L_2$.
Thus,   generally   $A_N$ converge in $L_2$ to $\int fd\mu \int gd\mu$.

Using induction, we  obtain the following statement.

\vspace{2mm}
\bf Theorem 1.3. \it Given  typical collection of  automorphisms
$T_1,\dots, T_m$,  the convergence    $$ \frac 1 N \sum_{n=1}^N  T^n_1f_1\dots T^n_mf_m \ \to_{L_2} \ \ \prod_{k=1}^m \int f_kd\mu$$
holds for  any functions $f_1,\dots, f_m\in L_\infty$. \rm

\vspace{2mm}
\bf Unknown, but typical examples. \rm The theory of generic automorphism properties 
is mainly of independent interest. Sometimes it shows the existence of automorphisms whose explicit construction is not easy to find. Let us give a modern example.
In connection with Kolmogorov's conjecture on the group property of   ergodic  automorphism spectrum, Oseledets in \cite{Os} proposed $\chi$-mixing: there is  a sequence $n_i$ for which the convergence 
$$ \mu(T^{n_i}A\cap B)\to \chi\mu(A)\mu(B) +(1-\chi)\mu(A\cap B)$$ 
holds for all measurable $A,B$. This property implies the disjointness of spectrum convolution powers. Stepin gave explicit examples and noted that $\chi$-mixing is  generic.
  Glasner, Thouvenot and Weiss \cite{GTW}   proved  that 
 the  generic automorphism  is relatively weakly mixing extension  of some of its nontrivial factors.      Examples of $\chi$-mixing skew products   with  relative weak mixing are    unknown, but typical. 

Further, in \S 2, \S 3 we consider metric spaces, which have not been sufficiently studied from the standpoint of the   generic property theory.
In \S 4, some unsolved problems on generic automorphisms are formulated.


\section{ The space ${\bf Aut}_\infty$ of infinite automorphisms}
Let $T$ be an automorphism of a space isomorphic to the line $\R$ with  Lebesgue measure. The theory of generic asymptotic properties of infinite automorphisms is basically similar to the theory of generic automorphisms of a probability space. 

Denote by $W(T)$ the weak closure of the group $\{T^n:\, n\in\Z\}$ in the algebra of bounded operators on $L_2$. Denote
$$Min(T) =\{\sum_{n\in\Z} a_n T^n, \ \ a_n\geq 0, \ \sum_{n\in\Z} a_n\leq 1\}.$$

An automorphism $T$ has \it rank one, \rm if there exists a sequence of partitions of the form
$$\{E_j, SE_j, S^2 E_{j},\dots, S^{h_j-1}E_j, \tilde E_j\},$$ tending to a partition into points.

We say that $T$ has \it $(n,n+1)$-approximation if     $E_j=E_j^1\uu E_j^2$, $\mu(E_j^1)=\mu (E_j^2)$, and $$ \mu(E_j^1\Delta T^{h_j}E_j^1) + \mu(E_j^2\Delta T^{h_j+1}E_j^2)\ \to\ 0.$$
 \rm
\vspace{2mm}
\bf Theorem 2.1. \it
The set $\left\{T\in {\bf Aut}_\infty\,:\, Min(T)\subset W(T)\right\}$ is generic. For any fixed sequence $n_i\to\infty$, the set of automorphisms of $T$ for which the sequence $T^{n_i}$ converges weakly has the first Baire category.

The generic automorphism $T\in {\bf Aut}_\infty$ has rank one,  cyclic $(n,n+1)$-approximation, its centralizer is the closure of the powers $T^n$,  Cartesian powers $T^{\times n}$ are ergodic,   spectra of symmetric tensor powers $T^{\odot n}$ are simple.  \rm

\vspace{2mm}
In \cite{20} there is an  automorphism $T$ with the following asymptotic property: for some sequence $h_j\to\infty$ for all sets $A$ of finite measure 
$$\mu(A\cap T^{h_j}A\cap T^{3h_j}A)\ \to {\mu(A)}/ 3    ,$$
$$\mu(A\cap T^{-h_j}A\cap T^{-3h_j}A)\ \to 0.$$
This  asymmetric property of the automorphism $T$ is generic, therefore, the asymmetry ($T$ is not isomorphic to $T^{-1}$) is generic as well.

\vspace{2mm}
\bf Theorem 2.2. \it The set of infinite automorphisms that are not isomorphic to their inverse is generic. \rm

\bf Density of conjugacy classes. \rm Theorems 2.1 and 2.2 are proved according to the classical scheme. We show that the set of automorphisms with the required invariant property is of type $G_\delta$. Then we find the required example $T$ and, importantly, show the density of its conjugacy class. For this, we use the following (infinite) analogue of the Rokhlin-Halmos lemma.

\vspace{2mm}
\bf Lemma 2.3. \it Let $T$ be an aperiodic infinite automorphism of the space
$(X,\B,\mu)$ with sigma-finite measure. 

(1) For every $n>0$ there exist
$N>n$ and measurable sets $B,B'$ such that
$$X=\left(\uu_{i=0}^{N-1}T^iB\right) \uu \left(\uu_{i=0}^{N}T^iB'\right),$$

(2) and   given $\eps>0$ and a set $A$ of finite measure, we have 
 $$\mu\left(A\cap (T^{N-1}B\cup T^{N-1}B'\cup T^NB')\right)<\eps.$$ \rm

\vspace{2mm}
Proof. The first part of the lemma is classical, it is established in the same way as in the probabilistic case: we find a high Kakutani tower, and then, using the mutual primality of the numbers $N$ and $N+1$, we obtain the desired partition. The second part of the lemma is deduced from the first as follows. We  take $N-1>3\mu(A)/\eps$ and, if $B,B'$ do not good, we consider $T^{-2n}B, T^{-2n}B'$ instead of $B,B'$. Obviously, for some $m$ with $0<2m<N-1$ we have
$$\mu\left(A\cap (T^{N-1-2m}B\cup T^{N-1-2m}B'\cup T^{N-2m}B')\right)<\eps.$$ Otherwise,  a contradiction: $mu(X\cap A)> \eps (N-1)/2 > 3 \mu(A)/2.$ Setting $B:=T^{-2m}B,$$B'= T^{-2m}B'$, we get (1),(2).

\vspace{2mm}
\bf Lemma 2.4. \it
Let $T\in {\bf Aut}_\infty$ and $S\in {\bf Aut}_\infty$ be aperiodic. For every $\delta>0$ there exists an automorphism of $R$ such that $\rho_\infty(T, RST^{-1})< \delta$.\rm

\vspace{2mm}
Proof. Let $T$ and $S$ be aperiodic. Let $A_1, A_2,\dots A_k$ be sets of finite measure that figure in the definition of the metric $\rho_\infty$. Denote their union by $A$, fix $\eps>0$.
Following Lemma 2.3, find the corresponding partitions for $T$ and $S$.
We choose a conjugation $R$ such that
$T$ and $RST^{-1}$ coincide on the set $$X\setminus (T^{N-1}B\cup T^{N-1}B'\cup T^NB').$$ Then, in the case of a sufficiently small $\eps$ we get  $\rho_\infty(T, RST^{-1})< \delta$.
In the general case, we  simultaneously approximate the aperiodic and periodic parts of the automorphism $T$.  We leave this as an exercise.

\section{Spaces $\bf {\bf Mix}$  of mixing automorphisms}
Define a metric $d_w$ on ${\bf Aut}(\mu)$:
$$ d_w(S,T)=\sum_{i,j=1}^\infty 2^{-i-j}\left|\mu(SA_i\cap A_j)-\mu(TA_i\cap A_j)\right|.$$
On the set of all mixing automorphisms of {\bf Mix} we define a metric $r$:
$$ r(S,T)= \rho(S,T) + \sup_{n>0}d_w(S^n,T^n).$$
In \cite{T07} it is shown that $({\bf Mix},r)$ is complete, separable metric space.

\vspace{3mm}
--- Singular spectrum, multiple mixing  are generic \cite{T06}.

---  Bashtanov proved \cite{Ba}    the genericity of rank 1, this implies  the triviality of the centralizer and the absence of factors. 

--- Simple spectra for all symmetric tensor powers are also generic (this follows from   \cite{Ba}, \cite{07}).

--- Generic mixing automorphisms  have   arbitrarily slow correlations \cite{24mz}.

--- Infinite Kushnirenko $P$-entropy is generic for mixing automorphisms \cite{21}.

\bf $P$-entropy. \rm We recall the definition of Kushnirenko $P$-entropy (see \cite{Ku},\cite{21}).  For a sequence $P$ of finite sets $P_j\subset \N$ and an automorphism $T$ of the probability space $(X,\mu)$, we define the entropy $h_P(T)$ as follows. Let $$h_j(T,\xi)=\frac 1 {|P_j|} H\left(\bigvee_{p\in P_j}T^p\xi\right),$$ where $\xi=\{C_1,C_2,\dots, C_n\}$ is a measurable partition of the set $X$. Recall that the partition entropy
is defined by the formula
$$ H(\xi)=-\sum_{i=1}^n \mu( C_i)\ln \mu( C_i).$$
Now we set
$$h_{P}(T,\xi)={\limsup_j} \ h_j(T,\xi),$$

$$h_{P}(T)=\sup_\xi h_{P}(T,\xi).$$
Note that in the case of $P_j=\{1,2,\dots, j\}$ the usual Kolmogorov entropy $h(T)$ coincides with $h_{P}(T)$.
In this article we will consider only sequences 
of expanding arithmetic progressions that are convenient for our purposes.
In \cite{21} it is shown that if $h(S)=0$, then $h_P(T)=0$ for some sequence $P_j=\{j,2j,\dots, L(j)j\}$, $L_j\to\infty.$

\

\vspace{2mm}
\bf Theorem 3.1. \it For any mixing automorphism $T$ of zero entropy, there exists a sequence $P$ of expanding arithmetic progressions such that $h_P(T)=0$ and for generic $S$ and any nontrivial partition $\xi$  we have $h_P(S,\xi)>0$.  \rm

\vspace{2mm}
This theorem has been proved in  \cite{21} for the space $\bf Aut$, but for $\bf Mix$ the  proof is the same.  Applying the arguments from the proof of Theorems 1.2,
taking into account the results on the disjointness of an automorphism with zero and completely positive entropy \cite{RT24}, we obtain the following  analogue of Theorem 1.3.

\vspace{2mm}
\bf Theorem 3.2. \it Given  typical in $\bf Mix$ collection of mixing automorphisms
$T_1,\dots, T_m$,  the convergence    $$ \frac 1 N \sum_{n=1}^N  T^n_1f_1\dots T^n_mf_m \ \to_{L_2} \ \ \prod_{k=1}^m \int f_kd\mu$$
holds for  any  $f_1,\dots, f_m\in L_\infty$. \rm

\vspace{2mm}

\vspace{3mm}
\section{  Space $ \bf Mix_\infty$ of infinite mixing automorphisms}
Let us consider  the standard  space with a sigma-finite measure $\mu_\infty$ and define on $ {\bf Mix}_\infty$ the following   metric $ r_\infty$:
$$ r_\infty(S,T)= \rho_\infty(S,T) + \sup_{n>0}d_w(S^n,T^n),$$
where
$$ d_w(S,T)=\sum_{i,j=1}^\infty 2^{-i-j}\left|\mu_\infty(SA_i\cap A_j)-\mu_\infty(TA_i\cap A_j)\right|,$$
and $\{A_i\}$ is a fixed sequence of sets,
dense in the family of all  finite measure sets.

\vspace{3mm}
\bf Theorem 4.1. (\cite{2407.21768}) \it The set of rank-one automorphisms is generic in the space $ {\bf Mix}_\infty$. \rm

\vspace{3mm}
The proof of this theorem uses the method of \cite{20}. However, instead of Ornstein's constructions, Sidon's constructions are now used (the definition is given below). By the result of \cite{RT}, rank-one for mixings implies another generic property: the triviality of the centralizer. The Kushnirenko entropy (for rapidly growing sequences) of Poisson suspensions over the generic automorphism in $ {\bf Mix}_\infty$ is infinite.

\vspace{3mm}
\bf Theorem 4.2. \it The set of asymmetric infinite mixings is generic. \rm

\vspace{3mm}
The proof uses the density of the conjugacy class of Sidon constructions (see below) and the following assertion. 

\vspace{3mm}
\bf Theorem 4.3 (\cite{2406.14390}). \it There exists a mixing Sidon automorphism $T$ such that for some sequences
$m(j,i)\, 1\leq i\leq r_j\to\infty$, for any set $A$ of finite measure the following convergences hold
$$\sum_{i=0}^{r_j-1}\mu(A\cap T^{m(j,i)}A\cap T^{-m(j,i+1)}A)\to \ \mu(A),$$
$$\sum_{i=0}^{r_j-1}\mu(A\cap T^{-m (j,i)}A\cap T^{m(j,i+1)}A)\to \ 0.$$
Such $T$ is not conjugate to $T^{-1}$. \rm

\vspace{2mm}
Modified Sidon constructions are used to prove the following theorems.

\vspace{2mm}
\bf Theorem 4.4. \it Let $0\neq f\in L_2(\mu)$ and $\psi(n)\to +0$ as
$n\to\infty$. The family of automorphisms $T\in {\bf Mix}_\infty$ for which the set $\{n\,:\, |(T^nf,f)|> \psi(n)\}$ is infinite contains a dense $G_\delta$-set. \rm

\vspace{2mm}
Similar result \cite{24mz} for the space ${\bf Mix}$ used Bashtanov's theorem on the density of the conjugacy class of any mixing $T\in {\bf Mix}$. 

\vspace{2mm}
\bf Theorem 4.5. \it For the generic infinite mixing automorphism $T$ its symmetric tensor powers have simple spectra.\rm

\vspace{2mm}
\bf Auxiliary constructions. \rm The above results on generic properties of elements of the space $ {\bf Mix}_\infty$ were obtained using Sidon constructions or their slight modifications. Recall their definition.
We fix a natural number $h_1$, a sequence $r_j\to\infty$ ($r_j$ is the number of columns into which the tower of stage $j$ is cut) and a sequence of integer vectors
$$ \bar s_j=(s_j(1), s_j(2),\dots, s_j(r_j-1),s_j(r_j)).$$
At step $j=1$ we have  a set of disjoint half-intervals
$E_1, SE_1,$ $\dots, S^{h_1-1}E_1 $. At step $j$, a system of non-intersecting half-intervals
$$E_j, SE_j, S^2 E_{j},\dots, S^{h_j-1}E_j$$ is defined
and on $ S^{n}E_j,$, excluding $ 1\leq n <h_j,$
transformation $S$ is a parallel translation. Such a set of half-intervals is called a tower of stage $j$, their union is denoted by $X_j$, and  also called a tower.

Stage $j+1$. Let us represent $E_j$ as a disjoint union of $r_j$ half-intervals
$E_j^1,E_j^2,\dots E_j^{r_j}$ of the same measure (length).
For each $i=1,2,\dots, r_j$ we consider the so-called column
$E_j^i, SE_j^i,\dots, S^{h_j-1}E_j^i.$
The union of these half-intervals is denoted by $X_{i,j}$.
To each column $X_{i,j}$ we add $s_j(i)$ non-intersecting half-intervals of the same measure as $E_j^i$, obtaining a set
$$E_j^i, SE_j^i, S^2 E_j^i,\dots, S^{h_j+s_j(i)-1}E_j^i$$
(all these sets are non-intersecting).
We denote $E_{j+1}= E^1_j$ and extend $S$ setting
$S^{h_j+s_j(i)}E_j^i = E_j^{i+1}$.Thus, the superstructure columns are built into a new tower of stage $j+1$, consisting of half-intervals
$$E_{j+1}, SE_{j+1}, S^2 E_{j+1},\dots, S^{h_{j+1}-1}E_{j+1},$$
where
$$ h_{j+1}=h_jr_j +\sum_{i=1}^{r_j}s_j(i).$$

The definition of the transformation $S$ at stage $j$ is preserved at all subsequent stages. As a result, on the space $X=\cup_j X_j$ we obtain an invertible transformation $S:X\to X$, preserving the standard Lebesgue measure on $X$.

\bf Sidon automorphisms. \rm Let the construction $S$ have the following property: \it the intersection
$X_j\cap S^mX_j$ for $h_{j}<m\leq h_{j+1}$ can be contained
only in one of the columns $X_{i,j}$ of the tower $X_j$. \rm Such transformations
are called \it Sidon. \rm The measure of space $X$ in this case is infinite. We can obtain a Sidon construction as follows: put $s_j(1)=4h_j,$ $s_j(i+1)=4s_j(i)$, $i<r_j\to\infty$.

Modifying a Sidon construction means changing
the Sidon parameters for some $i$, for example, for $(1-\eps_j)r_j<i<r_j$, $\eps_j\to +0$. The technique of modifications was used in \cite{07} for the mixing automorphosms of a probability space.
By analogy with this case,  the desired effect can be realized for an infinite automorphism without losing the density of its conjugacy class.

\section{Unsolved problems}
\bf 1. Homoclinic groups. \rm The weakly homoclinic group $WH(T)$ of an automorphism $T$ is defined as the set of automorphisms $S$ such that
$$ \left|\{n\,:\, \rho(Id\,,\, T^{-n}ST^n)<\eps, \ 1\leq n\leq N\}\right|/N\ \to \ 1, \ \eps>0, \  N\to\infty.$$  Ergodic Gaussian and Poisson suspensions have  ergodic weakly homoclinic groups \cite{19}.
 The group  $WH(T)$  is trivial for the generic $T\in \bf Mix$. It follows from the genericity of rank one \cite{Ba} and the fact that mixing rank one trasformations $T$ have trivial $WH(T)$ \cite{19}.  Infinite transformations have large ergodic weakly homoclinic groups.

Is the group $WH(T)$ trivial for   generic $ T\in {\bf Aut}$?

\vspace{2mm}
\bf 2. No prime factors. \rm
Is it true that any nontrivial factor of the generic automorphism has a proper factor too?

\vspace{2mm}
\bf 3. Factors with zero $P$-entropy. \rm Let  an automorphism for some $P$  have positive $P$-entropy and   some its nontrivial factor have  zero $P$-entropy. Is this situation generic?  

\vspace{2mm}
\bf 4. Orbits of compact sets. \rm Is the orbit $\{R^{-1}SR:\, R\in {\bf Aut}, S\in K\}$ of compact $K\subset{\bf Aut}$ a set of the first category? Perhaps automorphisms with infinitesimal positive entropy prevent a quick solution to this problem.
 The answer is positive, if  $K$ consists of automorphisms of zero entropy \cite{21} as well as  for $K$ consisting of positive entropy   automorphisms.  

\vspace{2mm}
\bf 5. Generic extensions. \rm Questions about generic properties of extensions of individual automorphisms (see \cite{GTW},\cite{23}) can be a source of nontrivial problems. 
Do Lebesgue  spectrum and multiple mixing be lifted  under the generic extension? \rm It is known that the singularity of spectrum and the mixing  are preserved \cite{23}.

\vspace{2mm}
\bf 6. Factors, roots. \rm Does  the generic infinite automorphism have nontrivial invariant sigma-algebras and  roots? Is it embedded in a flow?

\vspace{2mm}
\bf 7. Subtle spectral properties. \rm Under the assumption that  generic automorphism $T_1\in {\bf Aut_\infty}$ is included in a flow $\{T_t\}$ by the results of \cite{22}, the spectrum of all products of the form $T_{t_1}\otimes T_{t_2}\otimes\dots$ for $0<t_1<t_2<\dots$ will be simple.
Is it true that the Poisson suspension over such a flow inherits this property? Regardless of the inclusion of an automorphism in a flow, a similar question about the spectrum of a Gaussian suspension over it makes sense, since an ergodic Gaussian automorphism is always included in a continuum of flows.
The generic atomorphism $T\in {\bf Aut}$ is due to Lasaro  and de la Rue embedded in a flow \cite{LR}, does the discussed spectral property hold for this flow?

\vspace{2mm}
\bf 8. Analog of Bashtanov's theorem. \rm Is the conjugacy class of an infinite mixing automorphism  dense in $ {\bf Mix}_\infty$?

\large

\vspace{5mm}
Lomonosov Moscow State University

E-mail: vryzh@mail.ru

\end{document}